\documentstyle[11pt]{article}
\catcode`\@=11
\def\marginnote#1{}
\newcount\hour
\newcount\minute
\newtoks\amorpm
\hour=\time\divide\hour by60
\minute=\time{\multiply\hour by60 \global\advance\minute by-\hour}
\edef\standardtime{{\ifnum\hour<12 \global\amorpm={am}%
        \else\global\amorpm={pm}\advance\hour by-12 \fi
        \ifnum\hour=0 \hour=12 \fi
        \number\hour:\ifnum\minute<10 0\fi\number\minute\the\amorpm}}
\edef\militarytime{\number\hour:\ifnum\minute<10 0\fi\number\minute}

\def\draftlabel#1{{\@bsphack\if@filesw {\let\thepage\relax
   \xdef\@gtempa{\write\@auxout{\string
      \newlabel{#1}{{\@currentlabel}{\thepage}}}}}\@gtempa
   \if@nobreak \ifvmode\nobreak\fi\fi\fi\@esphack}
        \gdef\@eqnlabel{#1}}
\def\@eqnlabel{}
\def\@vacuum{}
\def\draftmarginnote#1{\marginpar{\raggedright\scriptsize\tt#1}}

\def\draft{\oddsidemargin -.5truein
        \def\@oddfoot{\sl preliminary draft \hfil
        \rm\thepage\hfil\sl\today\quad\militarytime}
        \let\@evenfoot\@oddfoot \overfullrule 3pt
        \let\label=\draftlabel
        \let\marginnote=\draftmarginnote
   \def\@eqnnum{(\theequation)\rlap{\kern\marginparsep\tt\@eqnlabel}%
\global\let\@eqnlabel\@vacuum}  }
\newcounter{app}
\newcounter{sapp}[app]
\def\theapp{\Alph{app}}
\newcommand{\app}[1]{
\refstepcounter{app}{\vspace{7mm}
\noindent\Large\bf Appendix
\theapp.
 \ #1 \par \vspace{5mm}}
\setcounter{equation}{0}
\def\theequation{\Alph{app}.\arabic{equation}}}

\def\a{\alpha}  \def\d{\delta} \def\o{\omega} \def\la{\lambda}
\def\b{\beta} \def\g{\gamma} \def\ot{\otimes}\def\gg{{\bf g}}

\newcommand{\Uqg}{\mbox{$U_{q}^{}(\widehat{{\bf g}} )\,$}}
\newcommand{\UqgD}{\mbox{$U_{q}^{(D)}(\widehat{{\bf g}} )\,$}}
\newcommand{\Uqdva}{\mbox{$U_{q}^{}(\widehat{{\bf sl}}_{{\bf 2}})\,$}}
\newcommand{\UqdvaD}{\mbox{$U_{q}^{(D)}(\widehat{{\bf sl}}_{{\bf 2}})\,$}}
\def\DYg{\widehat{DY}({\gg})^{(D)}}
\newcommand{\bn}{\begin{equation}}
\newcommand{\ed}{\end{equation}}
\newcommand{\rf}[1]{(\ref{#1})}
\newcommand{\W}{\mbox{$W_q(\widehat{{\bf g}})\,$}}
\newcommand{\Wg}{\mbox{$W_q(\widehat{{\bf g}})
\times U_{q}^{(D)}(\widehat{{\bf g}} )\,$}}
\newtheorem{proposition}{Proposition}
\renewcommand{\theequation}{{\thesection}.{\arabic{equation}}}

\def\hh{{\bf h}}
\def\R{{\cal R}}

\newcommand{\RR}{\mbox{$\overline{{\cal R}}$}}
\def\Ppexp{{\cal P}{\overrightarrow{\exp}}_{q^2}}
\def\Pqexp{{\cal P}{\overrightarrow{\exp}}_{q^{-2}}}
\def\Pppexp{{\cal P}{\overrightarrow{\exp}}_{q}}
\def\Pqqexp{{\cal P}{\overrightarrow{\exp}}_{q^{-1}}}
\def\Phexp{{\cal P}{\overrightarrow{\exp}}_{i\hbar}}
\def\K{{\cal K}}
\def\ZZ{{\bf Z}}
\def\CC{{\bf C}}
\def\he{{\hat e}}
\def\hf{{\hat f}}
\def\kk{{\hat a}}

\def\tih{\hbar}
\def\h{\hbar}
\def\bij{\frac{i\hbar(\a,\b)}{2}}
\def\ve{\varepsilon}
\def\Ayag{{\cal A}_{\tih}(\widehat{\bf g})^{(D)}}
\def\Ayadva{{\cal A}_{\tih}({\widehat{{\bf sl}}}_{{\bf 2}})^{(D)}}
\begin{document}
\begin{center}
\hfill ITEP-TH-23/98\\
\hfill RIMS-1198\\
\bigskip\bigskip
{\LARGE\bf Weyl group extension of quantized current algebras}\\
\vspace{0.3cm}
{\large\bf Jintai Ding$\ ^*\,$\footnote{E-mail: ding@kurims.kyoto-u.ac.jp}
and Sergei Khoroshkin$\ ^{\star *}$
\footnote{E-mail: khoroshkin@vitep5.itep.ru}}
\bigskip
\bigskip

$\quad^*${\it Research Institute for Mathematical Sciences, Kyoto University,
 Kyoto 606, Japan}

$\quad^\star${\it Institute of Theoretical \& Experimental Physics,
117259 Moscow, Russia}
\bigskip
\end{center}
\bigskip

\begin{abstract}
In this paper, we extend the Drinfeld current realization
of  quantum affine algebras  $U_q(\hat {\gg})$ and of the Yangians
in several directions:
we construct   current
operators for non-simple roots of
${\gg }$, define  a new braid group action in terms of the current operators
and  describe
the universal R-matrix for the corresponding
``Drinfeld'' comultiplication
in the form of infinite product and in the form of
certain integrals over current operators.
\end{abstract}
\setcounter{section}{-1}
\section{Introduction}

In the theory of simple Lie algebras, the theory of affine Lie algebras
$\hat {\gg}$
distinguishs  themselves from others due to the realization of the
affine Lie algebras in current operators (or loop realizations).
The current operator is an operator define as
$ x(z)= \Sigma _{n\in \ZZ} x_n z^{-n},$  where
 $x_n$ is an operator in the  algebra and $z$ is a formal variable.
One can associate to each root of $ {\gg}$  a current operator
$x_\alpha (z)$, which gives a complete description of
the whole  Lie algebra $\hat {\gg}$.
The advantage  of the current realization of affine Lie algebras
lies in  the fact that the formal variable $z$, which we
can also  treated as a number in
${\bf C}$,   allows us to   connect
the algebraic structures of affine Lie algebra with analytic,
geometric and topological structures, which is manifest best in the
theory of Knizhinik-Zamolochikov equations.
Current realization of  quantum affine algebras
was first given  by  Drinfeld  after a few years of the discovery of the
theory of quantum groups \cite{Drnew}. The connection of
the Drinfeld realization with the first realization of
quantum groups given by Jimbo and Drinfeld  and with other type
of realizations is  also highly nontrivial
\cite{Beck}\cite{KT4}\cite{KT5}\cite{DF}.
Drinfeld current realization
of $U_q(\hat {\gg})$  has since played extremely important
role in the theory of affine quantum algebras.

However from the structure point of view,
Drinfeld realization is not complete in the
sense that in this realization only
the current operators for the simple roots are given.
Also recently, it has gradually become clear that
there is a possibility to develop a theory of Drinfeld realization
 as a theory of quantized algebras of currents completely parallel to the
 theory of quantum groups.
  This is the problem we want to address in this paper.

The  receipt of contructing the new current operator follows from
a simple  idea in the theory of the affine Lie algebras,  which is to
use analytic properties of the matrix coefficients  of current operators.
Let $\alpha$, $\beta$ and $\alpha+\beta$ be positive roots of
${\gg}$. In the affine Lie algebra  $\hat {\gg}$, we have
$$ [X_\alpha(z), X_\beta(w)] = \delta(z/w) X_{\alpha+\beta}(z), $$
where $z$ and $w$ are formal variables and $\delta(z/w)= \Sigma_{n\in {\bf Z}}
 (z/w)^n$.  Let $V_\lambda$ be an highest weight representation of
$\hat {\gg}$, $V_\lambda^*$ its dual,  $v\in V_\lambda $ and
$ v^* \in V_\lambda^*$. We know that
$<v^*, X_\alpha(z)X_\beta(w)> $ as Laurent series
converges if we substitute $z$ and $w$ with
numbers in the region $|w| << |z|$.  From such a point of view,
we can see that  $X_{\alpha+\beta}(z)$ actually is an operator sitting
on the pole of $ X_\alpha(z)X_\beta(w)$, which is $(z-w)=0$. Thus analytically
we can rewrite it as
\bn
 \int_{ z\ {\rm around}\  w}  X_\alpha(z)X_\beta(w) dz/z =
 X_{\alpha+\beta}(w),     \label{I}
\ed
which is a well-known formula in field theory.
In the quantum case, we use a modification of \rf{I}.
The formula are  of the  following
type:
\bn
 Y_{\alpha+\beta}(w)=
\int_{ z\  {\rm around}\  wq^{\pm 1}}  Y_\alpha(z)Y_\beta(w) dz/z,
\label{II}
\ed
where $Y_\alpha(z)$ and  $Y_\beta(w)$ are current operators for the
simple roots from Drinfeld realization of $U_q(\hat {\gg})$.
The form of the commutation relations for the current operators of
Drinfeld realization and the properties of highest weight representations
give us the possibility to define precisely the analytic
continuation of the products of current operators.
This allows us to define R.H.S. of \rf{II}
as an element of algebra.
Applying this method further, we derive the whole set
of current operator corresponding to all the roots of ${\gg}$
for $U_q(\hat {\gg})$.

As we know in the classical case, the Weyl group is an
 important  and inseperable
structure in the theory of simple Lie algebra, once the
theory of root system is established. Going further in the develpment
 of the theory of quantized current algebras, we define an action of
 the braid group on \Uqg in terms of currents. This action uses,
 in particular, the composed root currents and does not coincide with
 Lusztig automorphisms \cite{Lusztig} of \Uqg. The automorphisms form a
 semidirect product of the braid group with the weight lattice.
  Beyond that we manage to  derive a new Hopf algebra
by combining both braid group automorphisms
and  $U_q(\hat {\gg})$.

 Current realization of quantum affine algebra is naturally supplied with
 its own quasitriangular Hopf algebra structure
(also given by Drinfeld \cite{Dr2}),
 which is different from that for quantum groups.
This  Hopf algebra
structure can be derived
as  a limit of  twists  of the Hopf algebra structure for
affine quantum groups given by Drinfeld and Jimbo
\cite{KT5}. The results of \cite{KT5}  also
give  the formula for the corresponding universal R-matrix ,
 which, however,  is
not given in the form of current operators and is not effictive for
 algebras of rank higher than two.
An attempt to derive the universal R-matrix for the
current realization  for
$U_q(\hat {\bf  sl}_2)$ was taken in \cite{Enriquez}.
  However,  their naive formula  is not correct.
 We derive a formula
 the universal R-matrix  using
a type of  noncommutative functional generalization of exponential
functions, where the choice of integration contours plays a fundamental
role. The formula is given first for the case of $U_q(\hat {\bf sl}_2)$
 (see \rf{En6}) and then is extended to general case by the Weyl group
 technique. Simulteniously we get a presentation of the universal $R$-matrix
 in a form of infinite product different from \cite{KT5} for generic $\gg$.
 For comparison, we give in Appendix the computation of quadratic terms
 of universal $R$-matrix for $U_q(\hat {\bf sl}_2)$ in both approaches.

The integral form of the universal $R$-matrix is given in a series of integrals
which, from one side, look as the complex-analytical counterpart of
 ordered exponentials over the real variable; the order of the arguments is
replaced by the
 precise choise of the poles inside the integration region. On the other hand,
 the properties of the integrals are direct generalizations of the properties
 of $q$-exponential functions \cite{FV}, \cite{KT1}, \cite{V}.
We observe them in a separate section.
 One can speculate further for  current generalizations of the quantum plane.

Our considerations are of quite general nature and can be applied as well
 to the Yangians and to the face type elliptic algebras. In section 5 we
 presents such results for the double of the Yangians. It is also clear
 how to extend them for non-simple laced case. Note also that due to the
 results in \cite{BLZ}, we can suppose that there should exist fumdamental
 presentations of the universal $R$-matrices in an integral
 form of ordered exponentials and the formulas which we present here could be
 viewed as  an approximation to these presentations.

\section{Composed root currents.}

Let $\gg$ be a simple Lie algebra of simple laced type and
\UqgD be the corresponding quantum affine algebra.
By \UqgD we mean that
 first, we define the algebra by  the ''new realization'' of \Uqg in
generating functions  by V.
 Drinfeld \cite{Drnew},
 second, this algebra has the Hopf algebra structure,
naturally attached to this
 realization \cite{Dr2} and third, we restrict this algebra to the
category of  highest weight representations.

 The generators of the algebra \UqgD are given by generating functions
 $e_{\a_i}(u)$, $f_{\a_i}(u)$ and $K^\pm_{\a_i}(u)$:
$$e_{\a_i}(z)=\sum_{k\in\ZZ}e_{\a_i,k}z^{-k},\qquad
f_{\a_i}(z)=\sum_{k\in\ZZ}f_{\a_i,k}z^{-k},$$
and
$$K^\pm_{\a_i}(z)=k_{\a_i}^{\pm 1}\exp\left(\pm(q-q^{-1})\sum_{n>0}
a_{i,\pm n}z^{\mp n}\right).$$
where $\a_i$, $i=1,\ldots,r$ are all simple roots of $\gg$.

 They satisfy the following defining relations
($\a,\b$ are simple roots of $\gg$):
\bn
(z-q^{(\a,\b)}w)e_\a(z)e_\b(w)=e_\b(w)e_\a(z)(q^{(\a,\b)}z-w)\ ,
\label{1}
\ed
\bn
(z-q^{-(\a,\b)}w)f_\a(z)f_\b(w)=f_\b(w)f_\a(z)(q^{-(\a,\b)}z-w)\ ,
\label{2}
\ed
\bn
\frac{(q^{c/2}z-q^{(\a,\b)}w)}{(q^{(\a,\b)+c/2}z-w)}K_\a^+(z)e_\b(w)=
e_\b(w)K_\a^+(z)\ ,
\label{3}
\ed
\bn
K_\a^-(z)e_\b(w)=\frac{(q^{(\a,\b)-c/2}z-w)}{(q^{-c/2}z-q^{(\a,\b)}w)}
e_\b(w)K_\a^-(z)\ ,
\label{4}
\ed
\bn
\frac{(q^{-c/2}z-q^{-(\a,\b)}w)}{(q^{-(\a,\b)-c/2}z-w)}K_\a^+(z)f_\b(w)=
f_\b(w)K_\a^+(z)\ ,
\label{5}
\ed
\bn
K_\a^-(z)f_\b(w)=\frac{(q^{-(\a,\b)+c/2}z-w)}{(q^{c/2}z-q^{-(\a,\b)}w)}
f_\b(w)K_\a^-(z)\ ,
\label{6}
\ed
\bn
\frac{(z-q^{(\a,\b)-c}w)(z-q^{-(\a,\b)+c}w)}
{(q^{(\a,\b)+c}z-w)(q^{-(\a,\b)-c}z-w)}
K_\a^+(z)K_\b^-(w)=K_\b^-(w)K_\a^+(z)
\label{7}
\ed
\bn
K_\a^\pm(z)K_\b^\pm(w)=K_\b^\pm(w)K_\a^\pm(z)
\label{7a}
\ed
\bn
[e_\a(z),f_\b(w)]=\frac{\delta_{\a,\b}}{q-q^{-1}}\left( \delta(z/q^cw)
K^+_\a(zq^{-c/2})-\delta(zq^c/w)K^-_\a(wq^{-c/2})\right)
\label{10}
\ed
and the Serre relations:
\bn
e_\a(w_1)e_\a(w_2)e_\b(z)-[2]_qe_\a(w_1)e_\b(z)e_\a(w_2)+
e_\b(z)e_\a(w_1)e_\a(w_2)+\ w_1\leftrightarrow w_2 =0
\label{8}
\ed
\bn
f_\a(w_1)f_\a(w_2)f_\b(z)-[2]_qf_\a(w_1)f_\b(z)f_\a(w_2)+
f_\b(z)f_\a(w_1)f_\a(w_2)
+\ w_1\leftrightarrow w_2=0
\label{9}
\ed \smallskip
 for any $\a,\b$, $(\a,\b)=-1$.
Here $\delta(z)=\sum_{k\in \ZZ}z^k$.

All the relations are the equations on formal power series, but for the
relations \rf{1}-\rf{7a} and \rf{8}, \rf{9},
 both sides of them are analytical functions with singularities at $z,w =0,
 \infty$ when they act on highest weight modules, while the relation \rf{10}
 is given in the sense of analytical continuation. In general we follow  the
 rule to understand the equalities without taking analytical continuations
 if the regions of Laurent expansions are not specified and there are
 no $\delta$- functions terms.

In a  highest weight  representation, the
 operator valued functions $e_\a(z)e_\b(w)$ are analytical in a region
 $|z|>|q^{(\a,\b)}w|$ with the products $e_{\a,n}e_{\b,m}$ being their Laurent
 coefficients in this region, and the relations \rf{1} allow us to define
 the analytical continuation of these functions into a region
 $|z|<|q^{(\a,\b)}w|$ as Laurent series
 \bn
\frac{q^{(\a,\b)}z-w}{z-q^{(\a,\b)}w}e_\b(w)e_\a(z).
\label{anal}
\ed
 The formula \rf{anal} can be considered as the definition of analytical
 continuations for the products of currents in algebra \UqgD.
 On the other hand, the restriction to the category of highest weight
 representations
 dictates, that we can consider the equalities in terms of normal ordered
 series of generators of the algebra like
 $$\sum_{k_1\leq k_2\leq\cdots \leq k_m}  C_{\a_i,k_i} e_{\a_1,k_1}
 e_{\a_2,k_2}\cdots e_{\a_m,k_m}.$$
The precise definition of the topological structure of the algebra \UqgD is
 given in \cite{KT5}, \cite{DK}.

This is the ideology we use throughout the paper.

For any two simple roots $\a$ and $\b$ let $e_{\a+\b}(z)$ and $f_{\a+\b}(z)$
 be the currents
\bn
e_{\a+\b}(z)\ =\  :e_\a(qz)e_\b(q^2z):\ ,
\qquad
f_{\a+\b}(z)\ =\  :f_\b(q^2z)f_\a(qz):\ ,
\label{12}
\ed
 that is
$$e_{\a+\b}(z)=\frac{1}{2\pi i}\oint\limits_{w\  {\rm around\ } qz}
e_\a(w)e_\b(q^2z)\frac{dw}{w},
\quad
f_{\a+\b}(z)=\frac{1}{2\pi i}\oint\limits_{w\  {\rm around\ } qz}
f_\b(q^2z)f_\a(w)\frac{dw}{w},$$
 where the integrands are taken in a sense of analytical continuations.
 Note that $e_{\a+\b}(z) \not= e_{\b+\a}(z)$.
 Performing the analytical continuation, we can rewrite
 the definition \rf{12} as
 a difference of two contour integrals:
\bn
e_{\a+\b}(z)=\frac{1}{2\pi i}\left(\oint_{C_\infty}e_\a(w)e_\b(q^2z)
\frac{dw}{w}-
\oint_{C_0}\frac{q^{-1}w-q^2z}{w-q^{}z}e_\b(q^2z)e_\a(w)
\frac{dw}{w}
\right)\ ,
\label{a}
\ed
\bn
f_{\a+\b}(z)=\frac{1}{2\pi i}\left(-\oint_{C_0}
f_\b(q^2z)f_\a(w)\frac{dw}{w}+
 \oint_{C_\infty}\frac{q^{-1}w-q^2z}{w-q^{}z}
f_\a(w)f_\b(q^2z)\frac{dw}{w}\right)\ ,
\label{b}
\ed
 where $C_0$ is a contour around zero such that the point $w=q^{}z$
is outside the
 contour while $C_\infty$ is a contour close to infinity, which includes zero
and the point  $w=q^{}z$. In other words,
\bn
e_\a(z)e_\b(w)- \frac{q^{-1}z-w}{z-q^{-1}w}e_\b(w)e_\a(z)=
\delta(zq/w)e_{\a+\b}(q^{-1}z)\ ,\quad |z|<|q^{-1}w|\ ,
\label{13}
\ed
\bn
f_\b(z)f_\a(w)- \frac{q^{}z-w}{z-q^{}w}f_\a(w)f_\b(z)=
\delta(z/qw)f_{\a+\b}(q^{-1}w)\ , \qquad |z|<|q^{}w|\ ,
\label{14}
\ed
In terms of Fourier modes
$e_{\g,n}=\frac{1}{2\pi i}\oint e_\g(w)w^{n}\frac{dw}{w}$,  the definitions
\rf{a} and \rf{b} can be rewritten as
\bn
e_{\a+\b,n}=q^{-2n-p}\left( e_{\a,p}e_{\b,n-p}-qe_{\b,n-p}e_{\a,p}-
(q-q^{-1})\sum_{k\geq 1}q^k e_{\b,n-p-k}e_{\a,p+k}\right)\ ,
\label{c}
\ed
\bn
f_{\a+\b,n}=q^{-2n-p}\left(-f_{\b,n-p} f_{\a,p}+q^{-1}f_{\a,p}f_{\b,n-p}-
(q-q^{-1})\sum_{k\geq 1}q^{-k} f_{\a,p-k}f_{\b,n-p+k}\right)
\label{d}
\ed
for any integer $p$.
\begin{proposition}
  The following relations are equivalent to the Serre relations:
\bn
e_\a(z)e_{\a+\b}(w)=\frac{qz-w}{z-qw}e_{\a+\b}(w)e_\a(z)\ ,
\quad |z|<|q^{}w|\ ,\label{15}
\ed
\bn
e_{\a+\b}(z)e_\b(w)=\frac{q^2z-q^{-2}w}{qz-q^{-1}w}e_\b(w)e_{\a+\b}(z)\ ,
\quad |z|<|q^{-2}w|\ ,\label{16}
\ed
\bn
f_\b(z)f_{\a+\b}(w)=\frac{q^{-2}z-q^2w}{q^{-1}z-qw}f_{\a+\b}(w)f_\a(z)\ ,
\quad |z|<|q^{2}w|\ ,\label{17}
\ed
\bn
f_{\a+\b}(z)f_\a(w)=\frac{q^{-1}z-w}{z-q^{-1}w}f_\a(w)f_{\a+\b}(z)\ ,
\quad |z|<|q^{-1}w|\ ,\label{18}
\ed
\end{proposition}
The equalities \rf{15}--\rf{18} are the equalities of Laurent series. In
 analytical language they are equalities of the functions, analytical
 in $\bigl({\bf C}^*\bigr)^2$ without restrictions to a region.

{\bf Proof.} Let us examine the relation \rf{17}. Due to the Serre
 relations and  \rf{2}, \rf{14}, we have
$$f_\a(w_1)f_\a(w_2)f_\b(z)-(q+q^{-1})f_\a(w_1)f_\b(z)f_\a(w_2)+
f_\b(z)f_\a(w_1)f_\a(w_2)=$$
$$\frac{(q-q^{-1})z(w_1-q^{-2}w_2)}{(w_2-q^{-1}z)(w_1-q^{-1}z)}
f_\a(w_1)f_\a(w_2)f_\b(z)+
f_{\a+\b}(q^{-1}w_1)f_\a(w_2)\d(z/qw_1)+$$
$$\frac{qz-w_1}{z-qw_1}f_\a(w_1)f_{\a+\b}(q^{-1}w_2)\d(z/qw_2)-
(q-q^{-1})f_\a(w_1)f_{\a+\b}(q^{-1}w_2)\d(z/qw_2)=$$
$$\frac{(q-q^{-1})z(w_1-q^{-2}w_2)}{(w_2-q^{-1}z)(w_1-q^{-1}z)}
f_\a(w_1)f_\a(w_2)f_\b(z)+$$
$$f_{\a+\b}(q^{-1}w_1)f_\a(w_2)\d(z/qw_1)-
\frac{q^{-1}z-q^2w_1}{z-qw_1}f_\a(w_1)f_{\a+\b}(q^{-1}w_2)\d(z/qw_2).$$
Then we do the same with arguments $w_1$ and $w_2$ interchanged and summ up.
The terms without delta functions dissapear while the coefficients before
 delta functions give the same equation
$$
f_{\a+\b}(q^{-1}w_1)f_\a(w_2)=
\frac{w_1-q^2w_2}{qw_1-qw_2}f_\a(w_2)f_{\a+\b}(q^{-1}w_1),
$$
which is equivalent to \rf{17}.
It is easy to see that the arguments can be reversed
 in order to deduce from \rf{17} one of the Serre relations.

  In particular \rf{15}--\rf{18} imply  an analog of PBW theorem for the
Fourier components of the currents for rank three algebra
$U_{q}^{(D)}(\widehat{{\bf sl}}_{{\bf 3}})\,$.

In the following we need also another pair of composed root currents:
\bn
\check{e}_{\a+\b}(z)\ =\  :e_\b(q^{-2}z)e_\a(q^{-1}z):\ ,
\qquad
\check{f}_{\a+\b}(z)\ =\  :f_\a(q^{-1}z)f_\b(q^{-2}z):\ ,
\label{11}
\ed
that is,
\bn
\check{e}_{\a+\b}(z)=\frac{1}{2\pi i}\left(-\oint_{{C'}_0}
e_\b(q^{-2}z)e_\a(w)
\frac{dw}{w}+
\oint_{{C'}_\infty}\frac{q^{}w-q^{-2}z}{w-q^{-1}z}e_\a(w)e_\b(q^{-2}z)
\frac{dw}{w}
\right)\ ,
\label{aa}
\ed
\bn
\check{f}_{\a+\b}(z)=\frac{1}{2\pi i}\left(\oint_{{C'}_\infty}
f_\a(w)f_\b(q^{-2}z)\frac{dw}{w}-
\oint_{{C'}_0}\frac{q^{}w-q^{-2}z}{w-q^{-1}z}
f_\b(q^{-2}z)f_\a(w)\frac{dw}{w}\right)\ ,
\label{bb}
\ed
where ${C'}_0$ is a contour around zero such that
the point $w=q^{-1}z$ is outside the
 contour while ${C'}_\infty$ is a contour close to infinity,
which includes zero
and the point  $w=q^{-1}z$.

In the following (see section 4), we define the composed root vectors for all
 roots of $\gg$. As usual, the definition depends on the choise of reduced
 decomposition for the longest element $w_0$ of the Weyl group of $\gg$.

\setcounter{equation}{0}
\section{The automorphisms}
 In this section we define a new braid group action on \UqgD using the
 current operators defined above.

{}For any simple root $\a$ let $T_\a$ be the following linear map:
\bn
T_\a e_\a(z)=f_\a(q^{-c}z){K_\a^+(q^{-c/2}z)}^{-1}\ ,
\label{19}
\ed
\bn
T_\a f_\a(z)={K_\a^-(q^{-c/2}z)}^{-1}e_\a(q^{-c}z)\ ,
\label{20}
\ed
\bn
T_\a K^\pm_\a(z)={K^\pm_\a(z)}^{-1}\ ,
\label{21}
\ed
and for the simple $\b$, $(\a,\b)=-1$
\bn
T_\a K^\pm_\b(z)=K^\pm_\a(q^{}z)K^\pm_\b(q^{ 2}z)\ ,
\label{22}
\ed
\bn
T_\a e_\b(z)=e_{\a+\b}(z)\ ,
\qquad
T_\a f_\b(z)=f_{\a+\b}(z)\ .
\label{24}
\ed
The rest of Drinfeld current operators for simple roots  is stable.
 We  introduce the maps
\bn
T_\a^{-1} e_\a(z)={K_\a^-(q^{c/2}z)}^{-1}f_\a(q^{c}z)\ ,
\label{19i}
\ed
\bn
T_\a^{-1} f_\a(z)=e_\a(q^{c}z){K_\a^+(q^{c/2}z)}^{-1}\ ,
\label{20i}
\ed
\bn
T_\a^{-1} K^\pm_\a(z)={K^\pm_\a(z)}^{-1}\ ,
\label{21i}
\ed
and for the simple $\b$, $(\a,\b)=-1$
\bn
T_\a^{-1} K^\pm_\b(z)=K^\pm_\a(q^{-1}z)K^\pm_\b(q^{- 2}z)\ ,
\label{22i}
\ed
\bn
T_\a^{-1} e_\b(z)=\check{e}_{\a+\b}(z)\ ,
\qquad
T_\a^{-1} f_\b(z)=\check{f}_{\a+\b}(z)\ .
\label{24i}
\ed

We also define the following automorphisms of affine shift.
For any simple root $\a$ we put
\bn
P_{\o_\a}e_\a(z)=ze_\a(z),\qquad
P_{\o_\a}f_\a(z)=z^{-1}f_\a(z) ,
\label{24a}
\ed
\bn
P_{\o_\a}K^\pm_\a(z)=q^{\mp c}K^\pm_\a(z).
\label{24b}
\ed
The rest of of Drinfeld current operators for simple roots  is stable.
Here $\o_\a$ denotes $\a$-th fundamental weight: $(\o_\a,\b) = \d_{\a,\b}$.

The maps $P_{\o_\a}$ and $P_{\o_\b}$ commute for different roots,
 so for any element $\la$ of a weight lattice $Q$ of Lie algebra $\gg$
 there is well defined automorphism $P_\la$.

\begin{proposition}\label{prop1}
The maps $T_\a^{\pm 1}$ and $P_\la$ can be extended to the automorphisms
 of the algebra \UqgD. They form a semidirect product of a
 braid group with a lattice $Q$, that is, $T_\a$  satisfy the braid group
relations
\bn
T_\a T_\b T_\a =T_\b T_\a T_\b,\qquad (\a,\b)=-1,
\label{25a}
\ed
\bn
T_\a T_\b=T_\b T_\a,\qquad (\a,\b)=0,
\label{25b}
\ed
\bn T_\a T_\a^{-1}=T_\a^{-1}T_\a =1
\ed
and reflect the shift automorphisms $P_\la$:
\bn
T_\a P_\la =P_{s_\a(\la)}T_\a\ .
\label{25c}
\ed
\end{proposition}
Here $s_\a$ is a reflection in $\hh^*$ with respect to a root $\a$:
 $s_\a(\la)=\la-2\frac{(\a,\la)}{(\a,\a)}\a$.

The proof of the automorphisms properties and of the braid group relations
 \rf{25a}, \rf{25b} and of \rf{25c} consists of direct verification for
 rank two and rank three algebras.

Let us make several remarks on the definition of automorphisms.
  One can see that the action of braid group generators $T_\a^{\pm 1}$
 is quite different from that of Lusztig automorphisms  \cite{Lusztig}.
 To the contrary,
 the affine shifts are actually present in usual quantum Weyl group of
 \Uqg , extended by the automorphisms of a Dynkin diagramm \cite{Beck}.
  There are also two subtle things in the formulas \rf{19}--\rf{24}
 and in the commutation relations of composed root currents,
 which do not appear in finite- dimensional case.
The first one is the shifts of the arguments in commutation
 relations \rf{16} and \rf{18} and in the definition of the action
 of braid group generators \rf{20}, \rf{20i} by  $q$. We have no
satisfactory explanation of this shift.

The second one is the inverse order in the relations \rf{15} and \rf{17},
 \rf{16} and \rf{18}. It is not as for usual Lusztig automorphisms which
 mean that the order of negative currents in PBW is opposite to the one
 for the positives.

 We denote by \W the corresponding group of automorphisms  and by
\Wg the associated semidirect product of an algebra
${\bf C}\cdot\W$ and of \UqgD. We define it by the relations
$$T_\a\cdot x\cdot T_\a^{-1}=T_\a(x),\qquad
P_\la\cdot x\cdot P_\la^{-1}=P_\la(x)$$
 where $x\in \UqgD$.
\setcounter{equation}{0}
\section{Universal $R$-matrices for \UqdvaD}
 There are two standard divisions of \UqdvaD into a sum of opposite
 Borel subalgebras. The two different comultiplication rules correspond
to them.
Let us first consider the possibility, which we denote further by (I),
 where the following comultiplication
rules take place (we write them for all \UqgD and omit an index of a simple
 root for \UqdvaD)
\bn
\Delta_{(I)} e_{\a}(z)=
e_{\a}(z)\ot 1 +K_{\a}^-(zq^{c_1/2})\ot e_{\a}(zq^{c_1}),
\label{Pa1}
\ed
\bn
\Delta_{(I)} f_{\a}(z)=
1\ot f_{\a}(z) +f_{\a}(zq^{c_2})\ot K_{\a}^+(zq^{c_2/2}) ,
\label{Pa2}
\ed
\bn
\Delta_{(I)} K_{\a}^+(z)=K_{\a}^+(zq^{c_2/2})\ot K_{\a}^+(zq^{-c_1/2}),
\label{Pa3}
\ed
\bn
\Delta_{(I)} K_{\a}^-(z)=K_{\a}^-(zq^{-c_2/2})\ot K_{\a}^-(zq^{c_1/2}),
\label{Pa4}
\ed
Let us treate
 the Hopf algebra \UqdvaD as a double of subalgebra, generated
 by $f(z)$ and $K^+(z)$.

The results by \cite{KT5}, see also \cite{KT4} imply that the universal
 $R$-matrix for the first case  with a restriction $|q|<1$ has a form
(in order to adapt the notations \cite{KT5} to \rf{1}--\rf{9}
 one should change $q$ to $q^{-1}$ everywhere in the formulas presented in
 \cite{KT5}):
\bn
\R_{(I)} =\K \cdot\RR^{}
\label{Pa16}
\ed
where
\bn
\K^{} =
q^{-\frac{h\ot h}{2}}q^{\frac{-c\ot d -d\ot c}{2}}
\exp\left(-(q-q^{-1})\sum_{n>0}\frac{n}{[2n]_q}a_n\ot a_{-n}\right)
q^{\frac{-c\ot d -d\ot c}{2}},\quad\ d=-z\frac{d}{dz}
\label{Pa17}
\ed
and
\bn
\RR^{} = \overrightarrow{\prod_{n\in\ZZ}}
\exp_{q^2}\left(-(q-q^{-1})f_{-n}\ot e_{n}\right).
\label{Pa18}
\ed
Here
$$\exp_p(x)=1+x+\frac{x^2}{1+p}+\frac{x^3}{(1+p)(1+p+p^2)}+\ldots
+\frac{x^n}{(n)_p!}+\ldots$$

 The   tensor of the Hopf pairing in \UqgD when
 it is considered as a double of of subalgebra, generated by $e(z)$ and
 $K^-(z)$ is, by definition,
$
\left({\R_{(I)}}^{21}\right)^{-1}.
$
 It has a form
\bn
\left({\R_{(I)}}^{21}\right)^{-1}=\RR'\cdot\K'
\label{Pa19}
\ed
where
\bn
\K'=\left(\K ^{21}\right)^{-1}=
q^{\frac{h\ot h}{2}}q^{\frac{c\ot d +d\ot c}{2}}
\exp\left((q-q^{-1})\sum_{n>0}\frac{n}{[2n]_q}a_{-n}\ot a_{n}\right)
q^{\frac{c\ot d +
d\ot c}{2}}
\label{Pa20}
\ed
and
\bn
\RR'=\left(\RR ^{21}\right)^{-1}=\overrightarrow{\prod_{n\in\ZZ}}
\exp_{q^{-2}}\left((q-q^{-1})e_{-n}\ot f_{n}\right).
\label{Pa21}
\ed

 Below we derive an integral form of the expressions
\rf{Pa18} and \rf{Pa21}, using the current operators $e(z)$ and $f(z)$.
 In order to do this we study the  pairing
 of two Hopf dual Borel subalgebras in \UqdvaD.

We will follow the convention that in the  pairing of Hopf algebras
$A$ and $A^0$, where $A^0$ is Hopf dual with opposite comultiplication,
 the elements of $A$ stand in the first place. In other words, we use the
 rules
\bn
<ab,x>=<a\ot b,\Delta'(x)>,\qquad <a, xy>=<\Delta(a),x\ot y>.
\label{Pa5}
\ed
Let us remind, that the commutation relations for the two subralgeras
 of a quantum double are given, due to Drinfeld, by the folowing rule:
\bn
 a\cdot b=<a^{(1)},b^{(1)}>\,<S^{-1}(a^{(3)}),b^{(3)}>\,
	b^{(2)}\cdot a^{(2)}
	\label{Pa6}
	\ed
	where
	$a\in A,\;b\in A^0,\; \Delta^2(x)=(\Delta\otimes Id)\Delta(x)=
x^{(1)}\otimes
	x^{(2)}\otimes x^{(3)}$, $S$ is  antipode in $A$.
Then, applying \rf{Pa6} to check of the relation
\rf{10} we have:
\bn
<f(z),e(w)>=-\frac{\delta(z/w)}{q-q^{-1}}
\label{Pa7}
\ed
Then the relation \rf{Pa5}, applied to commutation relation \rf{3}- \rf{6},
give the pairing for $K^\pm(u)$:
\bn
<K^+(z), K^-(w)>=\frac{q^2z-w}{z-q^2w}, \qquad |w|<|q^{-2}z|.
\label{Pa9}
\ed
 which means that
\bn
<k_\a, k_\a>=q^{-2},\quad <h,h>=-\frac{1}{\log q},\quad <a_n,a_{-n}>=
-\frac{[2n]_q}{n(q-q^{-1})}
\label{Pa13}
\ed
 where  $h=\log_q k_\a$.
Moreover, the rule \rf{Pa5} gives the following triangular decomposition
 of the pairing:
\bn
<K^+(z_1)f(z_2),e(w)>=<f(z_2),e(w)>,\quad <f(z), K^-(w_1)e(w_2)>=
<f(z), e(w_2)>.
\label{Pa15}
\ed
This produces  factorization \rf{Pa16} of the corresponding universal
$R$-matrix. The further study of the pairing of the products of fields
$K^+(z_i)$ and of $K^-(w_j)$ gives the factor $\K$. We omit
 these calculations.

By induction, one can deduce the following description for the pairing of
 the products of fields
$f(z_i)$ and $e(w_j)$:
$$<f(z_1)f(z_2)\cdots f(z_n),\: e(w_1)e(w_2)\cdots e(w_n)>=$$
\bn
\frac{1}{(q^{-1}-q^{})^n}\sum_{\sigma\in S_n}\prod_{i=1}^n\delta\Bigl(
\frac{z_i}{\sigma(w_j)}\Bigr)\prod_{i<j,\atop\sigma(i)>\sigma(j)}
g\Bigl(\frac{z_i}{z_j}\Bigr),
\label{Pa22}
\ed
where the function
\bn
g(z)=\frac{q^2-z}{1-q^2z}, \qquad |z|<q^{-2}
\label{Pa23}
\ed
is expanded in the mentioned region of analyticity $|z|<q^{-2}$. The
 relation \rf{Pa22} is an equality of formal power series.

The tensor $\RR^{} $ diagonilizes this pairing.
 It can be written in
the following ordered exponential form:
\begin{proposition}
\bn
\RR^{} = \Ppexp\left(\frac{(q^{-1}-q^{})}{2\pi i}
\oint f(z)\ot e(z)\frac{dz}{z}\right),
\label{En5}
\ed
where
$$\Ppexp\left(\frac{(q^{-1}-q)}{2\pi i}\oint f(z)\ot e(z)
\frac{dz}{z}\right)=1+$$
$$
 \sum_{n>0}\frac{(q^{-1}-q)^n}{n!(2\pi i)^n}\oint_{C_n}\frac{dz_n}{z_n}
\oint_{C_{n-1}}\frac{dz_{n-1}}{z_{n-1}}\cdots \oint_{C_1}\frac{dz_1}{z_1}
f(z_1)f(z_{2})\ldots f(z_n)\ot$$
\bn
 e(z_1)e(z_{2})\ldots e(z_n),
\label{En6}
\ed
 and
 the contours $C_k$ enclose an origin and are in a region where the points
 $z_i=q^{2}z_j$ are outside the contours $C_i$ and the points $z_i=q^{-2}z_j$
 are inside the contours $C_i$
\end{proposition}
One can choose the contours as follows. First, all the contours contain
 an origin. Then, the contour $C_1$ contains the points $q^{-2} z_j$, $j>1$
 and the points $q^{2} z_j$, $j>1$ are outside. The contour
$C_2$ contains the points $q^{-2} z_j$, $j>2$
 and the points $q^{2} z_j$, $j>2$ are outside and so on.

Note that in order to get an expression for $\RR $ in a component form,
 one should first analytically continuate the integrand of \rf{En6} into the
 region (suppose, for instance, $|q|>1$)
\bn
|q^{-2}z_j|<|z_i|<|q^{2}z_j|,\qquad i\not= j
\label{En6a}
\ed
and expand this continuation:
\bn
\prod_{i<j}\frac{q^{-2}-z_i/z_j}{1- q^{-2}z_i/z_j}f(z_1)f(z_{2})\ldots f(z_n)
\ot e(z_n)e(z_{n-1})\ldots e(z_1)
\label{En6b}
\ed
into Laurent series in the region \rf{En6a}. The result does not depend on a
restriction on $|q|$.

The second remark is about difference of the pictures for $|q|<1$ and
$|q|>1$. We see, that for $|q|>1$ the domain of integration has a natural form
 while for $|q|<1$ we have a sofisticated picture of the contours. On the other
 hand, the behaviour of the factors in the infinite products \rf{Pa17} and
 \rf{Pa18} is also different for different values of $|q|$.
The evaluation of the universal $R$-matrix \rf{Pa16}
in a tensor product of two-dimensional
 representations $V(z_1)\ot V(z_2)$ ($V$ has a basis $v_0,v_1$)  gives
  a triangular $R$-matrix with $\delta$-function term outside the
 diagonal \cite{KT5} for $|q|<1$ and is purely diagonal for $|q|>1$
 This is because the scalar factor
\bn
<v_0|\R^{(I) }| v_0>=\exp\left(-\sum_{n>0}
\frac{(q^n-q^{-n})^2}{q^{2n}-q^{-2n}}\frac{z^n}{n}
\right)
\label{Pa21a}
\ed
 coming from the action of bozons $a_n$ on $V(z)$ \cite {KT5}, \cite{Iohara}:
$$a_{\pm n}v_0=\frac{q^n-q^{-n}}{q-q^{-1}}\frac{z^{\pm n}}{n}v_0,\quad n>0$$
 admits different decompositions into infinite products for $|q|<1$ and
 $|q|>1$.

{\bf Proof.}
 The main step is the translation of formal series pairing \rf{Pa22} into
 analytical language. We can prove that the relation
$$<f(z_1)f(z_2)\cdots f(z_n),\: e(w_1)e(w_2)\cdots e(w_n)>=$$
\bn
\frac{1}{(q^{-1}-q^{})^n}\sum_{\sigma\in S_n}\prod_{i=1}^n\delta\Bigl(
\frac{z_i}{\sigma(w_j)}\Bigr)\prod_{i<j,\atop\sigma(i)>\sigma(j)}
\frac{q^{-2}-z_j/z_i}{1-q^{-2}z_j/z_i},
\label{En6c}
\ed
is valid as an equality of analytical functions in a region
 $|z_i|>|q^{-2}z_j|$,  $|w_i|<|q^{2}w_j|$, where $i<j$.
 In order to prove \rf{En6c}, we rewrite first the formal power series
 equality \rf{Pa22} in a form
$$<f(z_1)f(z_2)\cdots f(z_n),\: e(w_n)e(w_{n-1})\cdots e(w_1)>=$$
\bn
\frac{1}{(q^{-1}-q^{})^n}\sum_{\sigma\in S_n}\prod_{i=1}^n\delta\Bigl(
\frac{z_i}{\sigma(w_{n+1-j})}\Bigr)\prod_{i<j,\atop\sigma(i)>\sigma(j)}
g\Bigl(\frac{z_i}{z_j}\Bigr),
\label{Pa22d}
\ed
and then multiply both sides of \rf{Pa22d} by a formal power series
\bn
\prod_{i<j,\atop\sigma(i)>\sigma(j)}
\overline{g}\Bigl(\frac{w_i}{w_j}\Bigr),
\label{Pa22e}
\ed
where the function
\bn
\overline{g}(w)=\frac{q^{-2}-w}{1-q^{-2}w}, \qquad |w|<q^{2}
\label{Pa22g}
\ed
is considered in the  region  $|w|<q^{2}$.
Then in the rhs we have a cancellation of the rational functions over the same
 arguments, since they are expanded all in one direction, and, finally,
 $$<f(z_1)f(z_2)\cdots f(z_n),\: \prod_{i<j,\atop\sigma(i)>\sigma(j)}
\overline{g}\Bigl(\frac{w_i}{w_j}\Bigr)e(w_n)e(w_{n-1})\cdots e(w_1)>=$$
\bn
\frac{1}{(q^{-1}-q^{})^n}\sum_{\sigma\in S_n}\prod_{i=1}^n\delta\Bigl(
\frac{z_i}{\sigma(w_j)}\Bigr)\prod_{i<j,\atop\sigma(i)>\sigma(j)}
\frac{q^{-2}-z_j/z_i}{1-q^{-2}z_j/z_i},
\label{En22h}
\ed
 which is equivalent to \rf{En6c} since both sides converge in the region
$|z_i|>|q^{-2}z_j|$,  $|w_i|<|q^{2}w_j|$, where $i<j$.

The rest follows from technical calculations.
We should prove the following identities:
$$<\RR^{} ,\   e(w_1)e(w_2)\cdots e(w_n)\ot 1>=
 e(w_1)e(w_2)\cdots e(w_n)$$
and
$$<1\ot f(w_1)f(w_2)\cdots f(w_n),\ \RR^{} >=
 f(w_1)f(w_2)\cdots f(w_n).$$
Let us show,  how it works on quadratic terms.

 Taking, for instance, $e(w_1)e(w_2)\ot 1$ in the region
\bn
|w_1|<q^2|w_2|,
\label{En7}
\ed
  we have due to \rf{En6c}
$$\frac{(q-q^{-1})^2}{2(2\pi i)^2}<\oint_{C_2}\frac{dz_2}{z_2}
\oint_{C_1}\frac{dz_1}{z_1}
f(z_1)f(z_2)\ot e(z_1)e(z_2),\  e(w_1)e(w_2)\ot 1>=$$
$$
\frac{1}{2(2\pi i)^2}
 \oint_{C_2}\frac{dz_2}{z_2}
\oint_{C_1}\frac{dz_1}{z_1}\delta(z_1/w_1)\delta(z_2/w_2)
 e(z_1)e(z_2)+
$$
$$
\frac{1}{2(2\pi i)^2}\oint_{C_2}\frac{dz_2}{z_2}
\oint_{C_1}\frac{dz_1}{z_1}\frac{q^{-2}-z_1/z_2}{1- q^{-2}z_1/z_2}
 \delta(z_1/w_2)\delta(z_1/w_2) e(z_1)e(z_2)=
$$
$$\frac{1}{2}
 e(w_1)e(w_2) +\frac{q^{-2}-w_2/w_1}{2(1- q^{-2}w_2/w_1)}
e(w_2)e(w_1)$$
The second summand coincides with the first one due to the restriction
 \rf{En7} on the region of analyticity. This gives the first desired equality.
 For the second case, it is convenient to perform the pairing of the tensor
with $1\ot f(w_1)f(w_2)$ in
 a region $|w_1|>|q^{-2}w_2|$

A naive  version of \rf{En5} was suggested in \cite{Enriquez}.
 Note that one of the subtle points in the pairing \rf{Pa22} leading to a
 delicate answer is a
 restriction on a region in \rf{Pa23} which makes the pairing
 noncommutative.

The calculations for the case, where we consider \Uqdva as a double
 of subalgebra generated by the fields $e(z)$ and $K^-(z)$ are similar.
The pairing \rf{Pa22} modulo signs is given by the same formula, and the only
 (but very important) difference is that the coeffitients are in another
 regions of analyticity:
$$<e(w_1)e(w_2)\cdots e(w_n),\: f(z_1)f(z_2)\cdots f(z_n) >=$$
\bn
\frac{1}{(q-q^{-1})^n}\sum_{\sigma\in S_n}\prod_{i=1}^n\delta\Bigl(
\frac{z_i}{\sigma(w_j)}\Bigr)\prod_{i<j,\atop\sigma(i)>\sigma(j)}
{g}\Bigl(\frac{z_i}{z_j}\Bigr),
\label{En11}
\ed
where the function ${g}(z)=\frac{q^2-z}{1-q^2z}$
is considered in the  region  $|z|>q^{-2}$.

As a result, we have
\begin{proposition}
\bn
\RR'=\left(\RR ^{21}\right)^{-1}=
 \Pqexp\left(\frac{(q^{}-q^{-1})}{2\pi i}
\oint e(z)\ot f(z)\frac{dz}{z}\right),
\label{En13}
\ed
where
$$\Pqexp\left(\frac{(q^{}-q^{-1})}{2\pi i}\oint e(z)\ot f(z)
\frac{dz}{z}\right)=1+$$
$$
 \sum_{n>0}\frac{(q^{}-q^{-1})^n}{n!(2\pi i)^n}\oint_{C_n}\frac{dz_n}{z_n}
\oint_{C_{n-1}}\frac{dz_{n-1}}{z_{n-1}}\cdots \oint_{C_1}\frac{dz_1}{z_1}
e(z_1)e(z_{2})\ldots e(z_n)\ot$$
\bn
 f(z_1)f(z_{2})\ldots f(z_n),
\label{En14}
\ed
and
 the contours $C_k$ enclose the origin, such that the points
 $z_i=q^{-2}z_j$ are outside the contours $C_i$ and the points
 $z_i=q^{2}z_j$ are inside the contours $C_i$ in the same sense
as in Proposition 3.
\end{proposition}
The second comultiplication structure, which we again write for all
 \UqgD and drop an index of a simple root for \UqdvaD looks as follows:
\bn
\Delta_{(II)} e_{\a}(z)=
e_{\a}(z)\ot 1  +K_{\a}^+(zq^{-c_1/2})\ot e_{\a}(zq^{-c_1})  ,
\label{En17}
\ed
\bn
\Delta_{(II)} f_{\a}(z)=
1\ot f_{\a}(z) +e_{\a}(zq^{-c_2})\ot K_{\a}^-(zq^{-c_2/2}) ,
\label{En16}
\ed
\bn
\Delta_{(II)} K_{\a}^+(z)=K_{\a}^+(zq^{c_2/2})\ot K_{\a}^+(zq^{-c_1/2}),
\label{En19}
\ed
\bn
\Delta_{(II)} K_{\a}^-(z)=K_{\a}^-(zq^{-c_2/2})\ot K_{\a}^-(zq^{c_1/2}),
\label{En18}
\ed
 is connected to the first one by a simple twist over ''Cartan'' generators:
\bn
\Delta_{(II)}=\K ^{-1}{\Delta'}_{(I)}\K
\label{En18a}
\ed
 From \rf{En18a} it follows that
the universal $R$-matrix for the comultiplication (II) has a form:
\bn
\R_{(II)}=\RR^{21}\K,
\label{En20}
\ed
and
\bn
\left(\R_{(II)}^{21}\right)^{-1}=\left(\K^{21}\right)^{-1}\RR^{-1},
\label{En21}
\ed
where the entries of \rf{En20} and of \rf{En21} are given by \rf{Pa16},
 \rf{Pa17}, \rf{Pa20}, \rf{Pa21} and by \rf{En6}, \rf{En14}.

 \setcounter{equation}{0}
\section{Properties of ordered $q$-exponential integrals}
The main feature of the tensor $f(z)\ot e(z)$ which essentially implies
 the cocycle properties of the factor $\overline{\R}$ of the universal
$R$-matrix for \UqdvaD\, is that the function of two variables
 $G(z_1,z_2)=$
 $f(z_1)f(z_2)\ot e(z_1)e(z_2)$ is meromorphic in  $\bigl(\CC^*\bigr)^2$
 and has two simple poles at $z_1=q^{\pm 2}z_2$,
 the sum of whose residues is zero. Indeed, under assumption $|q|<1$
 the analytical continuation of the function
 $G(z_1,z_2)$ from the region $|z_1|>|q^{-2}z_2|$ into the region
 $|q^{2}z_2|<|z_1|<|q^{-2}z_2|$ is given as
$$G'(z_1,z_2)=
\frac{q^{-2}z_1-z_2}{z_1-q^{-2}z_2}f(z_2)f(z_1)\ot e(z_1)e(z_2),$$
and the analytical continuation of the function $G'(z_1,z_2)$ from the region
 $|q^{2}z_2|<|z_1|<|q^{-2}z_2|$ to the region $|z_1|<|q^{2}z_2|$ is
\bn
G''(z_1,z_2)=\frac{q^{-2}z_1-z_2}{z_1-q^{-2}z_2}\
\frac{q^{2}z_1-z_2}{z_1-q^{2}z_2}
f(z_2)f(z_1)\ot e(z_2)e(z_1) =f(z_2)f(z_1)\ot e(z_2)e(z_1).
\label{Gz}
\ed
From \rf{Gz} we see that integral of $G(z_1,z_2)$ over $z_1$ in a
  region $|z_1|>>|z_2|$ is the same as in a region $|z_1|<<|z_2|$, which
 means that the residues at the points $z_1=q^{\pm 2}z_2$ cancel.
Putting in a base the rule \rf{Gz} of analytical continuation, we formulate
 in this section an abstract definition of $q$-exponential integrals and
  prove their main algebraical properties.

Let $x(z)$ be operator valued function of complex variable $z$, satisfying
 the following conditions:

{\it (i)} For any integer $n$ the products $x(z_1)x(z_2)\cdots x(z_n)$ define
 an analytical function in a region $|z_1|>>|z_2|>>\ldots >>|z_n|>0$ which
 analytical continuation is a meromorphic function in $\bigl(\CC^*\bigr)^n$
 with only singularities being the simple
 poles on shifted diagonals $z_i=q^{}z_j$, $i\not= j$.

{\it (ii)} the analytical continuation of the product
$x(z_1)x(z_2)\cdots x(z_n)$ is a symmetric function in
$$X_n=\bigl(\CC^*\bigr)^n\setminus\bigcup_{i\not=j}\{ z_i=qz_j\}\ ,$$
in other words, the commutativity
\bn
x(z_1)x(z_2)=x(z_2)x(z_1)
\label{commun}
\ed
takes place in a sense of analytical continuations in $X_n$

For brevity we call such a function as almost commutative function with two
 poles. For any almost commutative function $x(z)$ with two poles we define
 ordered $q$-exponential integrals as follows:
\bn
\Pppexp\left(\frac{1}{2\pi i}\oint x(z)
\frac{dz}{z}\right)=1+
 \sum_{n>0}\frac{1}{n!(2\pi i)^n}\oint\limits_{C_1}\frac{dz_1}{z_1}
\cdots \oint\limits_{C_n}
\frac{dz_n}{z_n}\,
x(z_1)x(z_2)\cdots x(z_n),
\label{QE1}
\ed
\bn
\Pqqexp\left(\frac{1}{2\pi i}\oint x(z)
\frac{dz}{z}\right)=1+
 \sum_{n>0}\frac{1}{n!(2\pi i)^n}\oint\limits_{{C'}_1}\frac{dz_1}{z_1}
\cdots \oint\limits_{{C'}_n}
\frac{dz_n}{z_n}\,
x(z_1)x(z_2)\cdots x(z_n),
\label{QE2}
\ed
where the contours $C_i$ and ${C'}_{i}$ enclose the origin and for the
 contours $C_i$ the points $z_i=qz_j$ are outside the contour while the points
 $z_i=q^{-1}z_j$ are inside; for the
 contours ${C'}_i$ the points $z_i=qz_j$ are inside the contour while
the points  $z_i=q^{-1}z_j$ are outside. We do not discuss here the
 convergence of the rhs of \rf{QE1},
 \rf{QE2} and treate them as formal series.

We say also that two almost commutative functions $x_1(z)$ and $x_2(w)$ with
two
 poles {\it $q$-commute} and denote this as $[x(z),y(w)]_q=0$ if

{\it (iii)} For any
 integer $n$ and for any sequence $\iota=i_1,i_2,\ldots,i_n$, where $i_k=1,2$
 the products $x_{i_1}(z_{i_1,1})x_{i_2}(z_{i_2,2})\cdots x_{i_n}(z_{i_1,n})$
 are analytical functions in a region $|z_{i_1,1}|>>|z_{i_2,2}|>>\ldots >>
|z_{i_n,n}|>0$ which
 analytical continuation is a meromorphic function in $\bigl(\CC^*\bigr)^n$
 with only singularities being the simple
 poles on shifted diagonals $z_{1,k}=q^{}z_{1,l}$, $k\not= l$,
$z_{2,k}=q^{}z_{2,l}$, $k\not= l$ and $z_{1,k}=q^{}z_{2,l}$;

{\it (iv)}  the analytical continuation of the product
$x_{i_1}(z_{i_1,1})x_{i_2}(z_{i_2,2})\cdots x_{i_n}(z_{i_1,n})$
is a symmetric function in
$$X_\iota=\bigl(\CC^*\bigr)^n\setminus\left(\bigcup_{i,k,l:k\not=l}
\{ z_{i,k}=qz_{i,l}\}\cup_{k,l}\{ z_{1,k}=q^{}z_{2,l}\}\right) \ ,$$
where the action of the symmetric group is defined as
$$\sigma\ x_{i_1}(z_{i_1,1})\cdots x_{i_n}(z_{i_1,n})=
x_{i_{\sigma(1)}}(z_{i_{\sigma(1)},\sigma(1)})\cdots
x_{i_{\sigma(n)}}(z_{i_{\sigma(n)},\sigma(n)})\ .$$
In other words, the commutativity
\bn
x_1(z)x_2(w)=x_2(w)x_1(z)
\label{QE0}
\ed
takes place in a sense of analytical continuations in $X_\iota$.

One can  see that for any two almost commutative functions $x(z)$
 and $y(w)$ , which
 $q$-commute, their sum $x(z)+y(z)$ is also almost commutative with two poles,
 i.e., it satisfies the conditions {\it (i)} and {\it (ii)}. For two
 $q$-commuting functions $x(z)$
 and $y(w)$  we denote also by $:x(qz)y(z):$ the residue
$$:x(qz)y(z):=\frac{1}{2\pi i}\oint\limits_{w\ {\rm around}\ qz}
x(w)y(z)\frac{dw}{w}\ .$$
The following proposition is a straightforward corollary of the definitions
 above.
\begin{proposition}\label{addition}
{\rm (a)} For any almost commutative function $x(z)$ with two poles:
\bn
\left(\Pppexp\left(\frac{1}{2\pi i}\oint x(z)
\frac{dz}{z}\right)\right)^{-1}=
\Pqqexp\left(-\,\frac{1}{2\pi i}\oint x(z)
\frac{dz}{z}\right)
\label{QE3}
\ed
{\rm (b)} If two almost commutative functions $x(z)$ and $y(z)$ $q$-commute,
$[x(z),y(w)]_q=0$, then
\bn
\Pppexp\left(\frac{1}{2\pi i}\oint x(z)+y(z)
\frac{dz}{z}\right)=
\Pppexp\left(\frac{1}{2\pi i}\oint y(z)
\frac{dz}{z}\right)
\Pppexp\left(\frac{1}{2\pi i}\oint x(z)
\frac{dz}{z}\right)
\label{QE4}
\ed
{\rm (c)} For two almost commutative functions
$$
\Pppexp\left(\frac{1}{2\pi i}\oint x(z)
\frac{dz}{z}\right)
\Pppexp\left(\frac{1}{2\pi i}\oint y(z)
\frac{dz}{z}\right)=$$
\bn
\Pppexp\left(\frac{1}{2\pi i}\oint\Bigl( x(z)+y(z)+:x(qz)y(z):\Bigr)
\frac{dz}{z}\right)
\label{Faddeev}
\ed
\end{proposition}
 It is not difficult to see that the statements of Proposition 3 and of
Proposition 4 give another proof of the inversion property \rf{QE3} for
 almost commutative function $f(z)\ot e(z)$. The addition property
 \rf{QE4} is equivalent in this sense to the properties
$$(\Delta\ot id)\R=\R^{13}\R^{23}, \qquad
(id \ot\Delta)\R=\R^{13}\R^{12} $$
of the universal $R$-matrix for \UqdvaD and one could also consider \rf{QE4}
 as the corollary of Propositions 3 and 4.

There is also an analog of H'Adamard identity for ordered $q$-exponential
integrals. In order to write it down let us define the $q$-commutator operation
$$\left[\frac{1}{2\pi i}\oint x(z)
\frac{dz}{z},\ \frac{1}{(2\pi i)^n}\oint\limits_{C_{1}}\frac{dw_1}{w_1}
\cdots \oint\limits_{C_n}
\frac{dw_n}{w_n}\,
A(w_1,\ldots, w_n)\right]_q=$$
 $$=ad_q \frac{1}{2\pi i}\oint x(z)
\frac{dz}{z}\left(\frac{1}{(2\pi i)^n}\oint\limits_{C_{1}}\frac{dw_1}{w_1}
\cdots \oint\limits_{C_n}
\frac{dw_n}{w_n}\,
A(w_1,\ldots, w_n)\right)$$
as
$$\frac{1}{(2\pi i)^{n+1}}\left(\oint\limits_{C_{1}}\frac{dw_1}{w_1}
\cdots \oint\limits_{C_n}
\frac{dw_n}{w_n}\oint\limits_{C}\frac{dz}{z}
x(z)A(w_1,\ldots, w_n)-\right.
$$ $$\left. -
\oint\limits_{C_{1}}\frac{dw_1}{w_1}
\cdots \oint\limits_{C_n}
\frac{dw_n}{w_n}\oint\limits_{C'}\frac{dz}{z}
A(w_1,\ldots, w_n)x(z)\right),$$
where the contours $C$ and $C'$ enclose the origin and the points $z=qw_j$
are outside the contour $C$ and inside the contour $C'$ while the points
$z=q^{-1}w_j$ are inside  the contour $C$ and outside the contour $C'$.
\begin{proposition}
For any almost commutative function $x(z)$ and any operator $A$ the following
 identity of formal series takes place:
$$
\Pppexp\left(\frac{1}{2\pi i}\oint x(z)
\frac{dz}{z}\right)A\left(\Pppexp\left(\frac{1}{2\pi i}\oint x(z)
\frac{dz}{z}\right)\right)^{-1}=
A+\sum_{n>0}\frac{1}{n!}ad_q^n \frac{1}{2\pi i}\oint x(z)
\frac{dz}{z}(A).
$$
\end{proposition}
It is also proved by direct calculations where the idintities between the sums
 of $q$-numbers are replaced by the identities between different contours
integrals.
 \setcounter{equation}{0}
\section{ A Hopf structure of \Wg}
The semidirect product \Wg can be equiped with a structure of a Hopf algebra
 in the same manner as it was done for quntized envelopping algebras of
 Kac- Moody Lie algebras, see, e.g.\cite{Kirillov1}.

Let us recover an index of a simple root $\a$
 in previous section and denote the tensor $\RR $  given in
 \rf{Pa18}, \rf{En5}  as $\RR_{\a}$ and
 if it come from the fields $e_\a(u)$ and of
 $f_\a(u)$. We claim the following
\begin{proposition}
The semidirect product \Wg is a Hopf algebra with respect to
 the comultiplication structure which is:
$$\Delta(x)=\Delta_{(I)}(x)\qquad \mbox{for}\ x\in\UqgD,$$
$$\Delta(T_\a)=\left(T_\a\ot T_\a\right)\RR_{\a},\qquad
\Delta(T_\a^{-1})=\RR_{\a}^{-1}\left(T_\a^{-1}\ot T_\a^{-1}\right),$$
$$\Delta(P_\la)=P_\la\ot P_\la.$$
\end{proposition}

By means of the ordered exponentials of composed root currents we can
 construct the universal $R$-matrix for arbitrary \UqgD in simple laced case.
 Let us fix the reduced
 decomposition of the longest element $w_0$ of the Weyl group of
Lie algebra $\gg$ of the rank $r$
\bn
w_0=s_{\a_{i_1}}s_{\a_{i_{2}}}\cdots s_{\a_{i_N}}
\label{W3}
\ed
 and
 let $\check{e}_{\g_1}(u), \check{e}_{\g_{2}}(u), \ldots ,
\check{e}_{\g_N}(u)$
 and
$\check{f}_{\g_1}(u), \check{f}_{\g_{2}}(u), \ldots ,\check{f}_{\g_N}(u)$ be
 the composed root currents given by the prescription:
$$\check{e}_{\g_N}(u)=e_{\a_{i_N}}(u),\qquad
  \check{f}_{\g_N}(u)=f_{\a_{i_n}}(u),$$
$$\check{e}_{\g_{N-1}}(u)=T_{\a_{i_N}}^{-1}e_{\a_{i_{N-1}}}(u),\quad
  \check{f}_{\g_{N-1}}(u)=T_{\a_{i_N}}^{-1}f_{\a_{i_{N-1}}}(u)$$
$$\check{e}_{\g_{N-2}}(u)=T_{\a_{i_N}}^{-1}T_{\a_{i_{N-1}}}^{-1}
e_{\a_{i_{N-2}}}(u),\quad
\check{f}_{\g_{N-2}}(u)=T_{\a_{i_N}}^{-1}T_{\a_{i_{N-1}}}^{-1}
f_{\a_{i_{N-2}}}(u)$$
and so on.
 Let also $\RR_{\g_i}$ be a tensor $\RR $ constructed from the composed fields
 $\check{e}_{\g_i}(u)$ and $\check{f}_{\g_i}(u)$ and
\bn
\K^{} =
q^{-t}q^{\frac{-c\ot d -d\ot c}{2}}
\prod_{n>0}\exp\left(-n(q-q^{-1})\sum_{i,j=1}^r d_{i,j}^{(n)}
a_{i,n}\ot a_{j,-n}\right)
q^{\frac{-c\ot d -d\ot c}{2}}
\label{W4}
\ed
where $t=\sum h_i\ot h^i$ is an invariant tensor in tensor square of
 Cartan subalgebra $\hh\ot\hh$ and $d_{i,j}^{(n)}$ is an inverse matrix to
 $$b_{i,j}^{(n)}=[n(\a_i,\a_j)]_q.$$
We claim that
\begin{proposition}
The tensor
\bn
\R=\K \overrightarrow{\prod}_{1\leq i\leq N}\RR_{\g_i}
\label{W6}
\ed
is the universal $R$-matrix for the algebra \UqgD with comultiplication
$\Delta_{(I)}$.
\end{proposition}
The logic of the proofs of Propositions 5 and 6 is standard..
  The role of the tensor
$\K $ from \rf{W4}
is the same as of $q^t$ in a Kac-Moody case: one can check directly,
 that it satisfies the condition \rf{En18a}.
 This condition reduces the proof of Propositions
 to the check of certain cocycle
 conditions on ordered $q$-exponential integrals. The new input is
 algebraical properties of $q$-exponential integrals which are the
generalizations of the properties of usual $q$-exponential functions,
 (see, e.g., \cite{KT1}) and are equivalent to the cocycle conditions.
  Note also, that in the formula
 \rf{W6} we can freely choose the presentation of the factors $\RR_{\g_i}$
 either in a form of infinite product or in an integral form. The integral
form is  more convenient for the proof of propositions. Still, the infinite
 product variant of \rf{W6} gives a presentation of the universal $R$-matrix
 for \UqgD, different from \cite{KT5} if rank $\gg$ greater then one.


\setcounter{equation}{0}
\section{The Yangian version }
It is not difficult to rewrite the whole story for the double of the Yangian.
Let us choose here for definiteness the presentaion of the double of the
 Yangian $\Ayag$ in Fourier integrals \cite{KLP5} where its
generating functions
 are given by the prescriptions
$$e_{\a_i}(u)=\int_{-\infty}^{+\infty}d\la\ e^{-i\la u}\he_{\a_i,\la}\ ,\qquad
 f_{\a_i}(u)=\int_{-\infty}^{+\infty}d\la\ e^{-i\la u}\hf_{\a_i,\la}\ ,$$
 $$K^\pm_{\a_i}(u)=\exp \h\int^{\pm\infty}_0d\la\
e^{-i\la u}\kk_{i,\la}\ .$$
 and satisfy the additive version of the relations \rf{1}--\rf{9}:
$$
 (u-v+i\tih(\a,\b)/2)e_\a(u)e_\b(v)=
  e_\b(v)e_\a(u)(u-v-i\tih(\a,\b)/2)\ ,
$$ $$
(u-v-i\tih(\a,\b)/2)f_\a(u)f_\b(v)=
 f_\b(v)f_\a(u)(u-v+i\tih(\a,\b)/2)\ ,
$$ $$
{(u-v+\bij-\frac{i\tih c}{4})
\over (u-v-\bij- \frac{i\tih c}{4})}K^+_\a(u)e_\b(v)=
e_\b(v)K^+_\a(u)\ ,
$$ $$
K^-_\a(u)e_\b(v)=
{ (u-v-\bij+\frac{i\tih c}{4})
\over
(u-v+\bij+\frac{i\tih c}{4})
}\
e_\b(v)K^-_\a(u)\ ,
$$ $$
{(u-v-\bij+\frac{i\tih c}{4})
\over(u-v+\bij+\frac{i\tih c}{4})}K^+_\a(u)f_\b(v)=f_\b(v) K^+_\a(u)\ ,
$$ $$
K^-_\a(u)f_\b(v)=
{ (u-v+\bij- \frac{i\tih c}{4})
\over
(u-v-\bij- \frac{i\tih c}{4})
}\
f_\b(v) K^-_\a(u)\ ,
$$ $$
\frac{(u-v+\bij-\frac{i\tih c}{2})(u-v-\bij+\frac{i\tih c}{2})}
{(u-v+\bij+\frac{i\tih c}{2})(u-v-\bij-\frac{i\tih c}{2})}
K_\a^+(u)K_\b^-(v)=K_\b^-(v)K_\a^+(u),
$$ $$
K_\a^\pm(u)K_\b^\pm(v)=K_\b^\pm(v)K_\a^\pm(u)
 $$ $$
{[}e_\a(u),f_\b(v){]}= \frac{i\delta_{\a,\b}}{\tih}\left[
\delta\left(u-v+\frac{ic\tih}{2}\right)
K^+_\a\left(u+\frac{ic\tih}{4}\right)-
\delta\left(u-v-\frac{ic\tih}{2}\right)
K^-_\a\left(v+\frac{ic\tih}{4}\right)
\right]  ,
$$
and the Serre relations
$$
e_\a(u_1)e_\a(u_2)e_\b(v)-2e_\a(u_1)e_\b(v)e_\a(u_2)+
e_\b(v)e_\a(u_1)e_\a(u_2)+\ u_1\leftrightarrow u_2 =0
$$ $$
f_\a(u_1)f_\a(u_2)f_\b(v)-2f_\a(u_1)f_\b(v)f_\a(u_2)+
f_\b(v)f_\a(u_1)f_\a(u_2)+\ u_1\leftrightarrow u_2 =0
$$ \smallskip
 for any $\a,\b$, $(\a,\b)=-1$.
Here $\delta$ function is defined as
\bn
 \delta(u-v)=\lim_{\epsilon\to0} \left[ {1\over
u-v-i\epsilon}-{1\over u-v+i\epsilon}\right]= i \int\limits_{-\infty}^{\infty}
 d\la\ {\rm e}^{-i\la
(u-v)}\ .
\label{delta-con}
\ed
The highest weight representations are characterized by the property that
 $$\he_{\a,\la}v=\hf_{\a,\la}v=\kk_{\a,\la}v=0$$
for any vector $v$ in the representation space for $\la$ big enough.
 It means, in particular, that
 operator valued functions $e_\a(u)e_\b(v)$ are analytical in a region
 ${\rm Im}\, u>>{\rm Im}\, v$ , and the relations \rf{1} allow to define
 the analytical continuation of these functions into a region
 ${\rm Im}\, u<{\rm Im}\,( v-i\tih(\a,\b)/2)$ as
 $$
e_\a(u)e_\b(v)=
\frac{(u-v-\bij)}{(u-v+\bij)}  e_\b(v)e_\a(u).
 $$
The definition of composed root vectors has a form:
$$
e_{\a+\b}(v)\ =\  :e_\a(v-i\tih/2)e_\b(v-i\tih):\ ,
\qquad
f_{\a+\b}(v)\ =\  :f_\b(v-i\tih)f_\a(v-i\tih/2):\ ,
$$
$$
\check{e}_{\a+\b}(v)\ =\  :e_\b(v+i\tih)e_\a(v+i\tih/2):\ ,
\qquad
\check{f}_{\a+\b}(v)\ =\  :f_\a(v+i\tih/2)f_\b(v+i\tih):\ ,
$$
 that is
\bn
e_{\a+\b}(v)=\frac{1}{2\pi i}\left(\int_{C_{+i\infty}}e_\a(u)e_\b(v-i\tih)
{du}-
\int_{C_{-i\infty}}\frac{(u-v+\frac{3i\tih}{2})}{(u-v+\frac{i\tih}{2})}
e_\b(v-i\tih)e_\a(u)
{du}
\right)\ ,
\label{Y10}
\ed
\bn
f_{\a+\b}(v)=\frac{1}{2\pi i}\left(-\int_{C_{-i\infty}}
f_\b(v-i\tih)f_\a(u){du}+
 \int_{C_{+i\infty}}\frac{(u-v+\frac{3i\tih}{2})}{(u-v+\frac{i\tih}{2})}
f_\a(u)f_\b(v-i\tih){du}\right)\ ,
\label{Y11}
\ed
 where $C_{+i\infty}$ is a contour going from $-\infty$ to $+\infty$ in such
 a way that the points $u=v-i\tih/2$ and $u=+i\infty$ are on the different
 sides of the contour and $C_{-i\infty}$ is a contour going from $-\infty$ to
$+\infty$ in such
 a way that the points $u=v-i\tih/2$ and $u=-i\infty$ are on the different
 sides of the contour;
$$
\check{e}_{\a+\b}(v)=\frac{1}{2\pi i}\left(-\int_{{C'}_{-i\infty}}e_\b(v+i\tih)
e_\a(u)
{du}+
\int_{{C'}_{+i\infty}}\frac{(u-v-\frac{3i\tih}{2})}{(u-v-\frac{i\tih}{2})}
e_\a(u)e_\b(v+i\tih)
{du}
\right)\ ,
$$ $$
\check{f}_{\a+\b}(v)=\frac{1}{2\pi i}\left(\int_{{C'}_{+i\infty}}
f_\a(u)f_\b(v+i\tih){du}-
 \int_{{C'}_{-i\infty}}\frac{(u-v-\frac{3i\tih}{2})}{(u-v-\frac{i\tih}{2})}
f_\b(v+i\tih)f_\a(u){du}\right)\ ,
$$
 where ${C'}_{+i\infty}$ is a contour going from $+\infty$ to $-\infty$ in such
 a way that the points $u=v+i\tih/2$ and $u=+i\infty$ are on the different
 sides of the contour and ${C'}_{-i\infty}$ is a contour going from $+\infty$
to $-\infty$ in such
 a way that the points $u=v+i\tih/2$ and $u=-i\infty$ are on the different
 sides of the contour.

The Serre relations are equivalent to
$$
e_\a(u)e_{\a+\b}(v)=\frac{u-v-i\tih/2}{u-v+i\tih/2}e_{\a+\b}(v)e_\a(u)\ ,
\qquad {\rm Im\; }u<\ {\rm Im\; }(v-i\tih/2)\ ,
$$ $$
e_{\a+\b}(u)e_\b(v)=\frac{u-v-2i\tih}{u-v-i\tih}e_\b(v)e_{\a+\b}(u)\ ,
\qquad {\rm Im\; }u<\ {\rm Im\; }(v+i\tih)\
$$ $$
f_\b(u)f_{\a+\b}(v)=\frac{u-v+2i\tih}{u-v+\tih}f_{\a+\b}(v)f_\a(u)\ ,
\qquad {\rm Im\; }u<\ {\rm Im\; }(v-i\tih)\
$$ $$
f_{\a+\b}(u)f_\a(v)=\frac{u-v+i\tih/2}{u-v-i\tih/2}f_\a(v)f_{\a+\b}(u)\ ,
\qquad {\rm Im\; }u<\ {\rm Im\; }(v+i\tih/2)\
$$
and the Weyl group automorphisms have a form
$$
T_\a e_\a(u)=f_\a(u+\frac{i\tih c}{2}){K_\a^+(u+\frac{i\tih c}{4})}^{-1}\ ,
$$ $$
T_\a f_\a(u)={K_\a^-(u+\frac{i\tih c}{4})}^{-1}e_\a(u+\frac{i\tih c}{2})\ ,
$$ $$
T_\a K^\pm_\a(u)={K^\pm_\a(u)}^{-1}\ ,\qquad
T_\a K^\pm_\b(u)=K^\pm_\a(u-\frac{i\tih }{2})K^\pm_\b(u-i\tih) ,
\quad (\a,\b)=-1,
$$ $$
T_\a e_\b(u)=e_{\a+\b}(u)\ ,
\qquad
T_\a f_\b(u)=f_{\a+\b}(u)\ . \quad (\a,\b)=-1,
$$ $$
T_\a^{-1} e_\a(u)={K_\a^-(u-\frac{i\tih c}{4})}^{-1}
f_\a(u-\frac{i\tih c}{2})\ ,
$$ $$
T_\a^{-1} f_\a(u)=e_\a(u-\frac{i\tih c}{2})
{K_\a^+(u-\frac{i\tih c}{4})}^{-1}\ ,
$$ $$
T_\a^{-1} K^\pm_\a(u)={K^\pm_\a(u)}^{-1}\ ,\quad
T_\a^{-1} K^\pm_\b(u)=K^\pm_\a(u+\frac{i\tih }{2})K^\pm_\b(u+i\tih)\ ,
\quad (\a,\b)=-1,
$$ $$
T_\a^{-1} e_\b(u)=\check{e}_{\a+\b}(u)\ ,
\qquad
T_\a^{-1} f_\b(u)=\check{f}_{\a+\b}(u)\ \quad (\a,\b)=-1,
$$
\bn
P_{\ve\o_\a}e_\a(u)=e^{i\ve u}e_\a(u),\qquad
P_{\ve\o_\a}f_\a(u)=e^{-i\ve u}f_\a(u) , \quad \ve\in {\bf R}
\label{Y24a}
\ed
\bn
P_{\ve\o_\a}K^\pm_\a(u)=e^{\pm\ve \tih/2}K^\pm_\a(u).
\label{Y24b}
\ed
 Here the weight lattice $Q$ is replaced by its real form $Q_{{\bf R}}=
 Q\ot_{{\bf Z}}{{\bf R}}$.
All the statements of the first  sections remain valid in the Yangian
 case.

The universal $R$-matrix for $\Ayadva$ with the comultiplication rule
$$
\Delta e_{\a}(u)= e_{\a}(u)\ot 1 +K_{\a}^-(u-ic_1\tih/4)\ot
e_{\a}(u-ic_1\tih/4),
$$ $$
\Delta f_{\a}(u)= 1\ot f_{\a}(u) +
f_{\a}(u-ic_2\tih/4)\ot K_{\a}^+(u-ic_2\tih/4) ,
$$ $$
\Delta  K_{\a}^+(u)=K_{\a}^+(u-ic_2\tih/4)\ot K_{\a}^+(u+ic_1\tih/4),
$$ $$
\Delta  K_{\a}^-(u)=K_{\a}^-(u+ic_2\tih/4)\ot K_{\a}^-(u-ic_1\tih/4)
$$
can be presented as
$$
\R=\K\RR
$$
where $(d=\frac{d}{du})$
\bn
\K= e^{i\h( d\ot c+c\ot  d)/4}\exp \left(-\int_0^{+\infty} d\la\  {\h^2\la
\over 2\,{\rm sh}\,\h\la}
\kk_{\la}\ot \kk_{-\la}\right)\ e^{i\h( d\ot c+c\ot  d)/4.}
\label{R00}
\ed
\bn
\RR =\  \overrightarrow{P} \exp
\left(-\h\int_{-\infty}^{+\infty} d\la\,
\hf_{-\la}\ot \he_{\la}\right)\ .
\label{Rpm}
\ed
or, equivalently,
\bn
\RR^{} = \Phexp\left(\frac{\tih}{2\pi }
\int_{-\infty}^{+\infty} f(u)\ot e(u){du}\right),
\label{Yn5}
\ed
where
$$\Phexp\left(\frac{\tih}{2\pi }\int_{-\infty}^{+\infty} f(u)\ot e(u)
{du}\right)=1+$$
$$
 \sum_{n>0}\frac{\tih^n}{n!(2\pi )^n}\int_{C_n}{du_n}
\int_{C_{n-1}}{du_{n-1}}\cdots \oint_{C_1}{du_1}
f(u_1)f(u_{2})\ldots f(u_n)\ot$$
\bn
 e(u_1)e(u_{2})\ldots e(u_n),
\label{Yn6}
\ed
 and
 the contours $C_k$ are going from $-\infty$ to $+\infty$ in such a way that
 the points $u_k=u_j+i\tih$ and the point $-i\infty$ are on on the one side
of the contour $C_k$ while the points  $u_k=u_j-i\tih$ and the point $+i\infty$
 are on the other side of the contour.
 The arguments of the section 4 can be repeated here and they give the
 universal $R$-matrix for $\Ayag$ for any simple laced $\gg$.

{}For the usual presentation of the Yangian double $\DYg$ by Laurent series
  (see, e.g., \cite{KLP5}), we use instead of the marked points $\pm\infty$
zero and infinity,
 as in the case of quantum affine algebra. The formulas are the same except
 \rf{Y24a} and \rf{Y24b} where the lattice appear again and the definition of
 affine shifts is more close to the case of quantum affine algebras.
 In an analogous manner we can treate the face type elliptic algebras and
 the elliptic algebras in scaling limit with Drinfeld's type comultiplication.
 We will observe it a separate note.

\bigskip
{\Large\bf Acknowlegements}
\medskip

This work was done when the second author visited RIMS.
He use the opportunity to appreciate the Institute and Prof. T. Miwa for the
kind hospitality.
The authors are grateful to Profs. V. Bazhanov, L.D. Faddeev, B. Feigin,
S. Kharchev, S. Pakuliak and F. Smirnov for the discussions.
We thak also Profs. S.Kharchev and
 S.Pakuliak for the help  in checking the key formulas. S.Kh. was supported
by INTAS grant 93-10183, RFBR grant 98-01-00303 and
grant  96-15-96455  for support of scientific schools.

\setcounter{equation}{0}
\app{ The braid group relation}
Let us examine the identity
\bn
T_\a T_\b T_\a\bigl(e_\b(z)\bigr)=T_\b T_\a T_\b\bigl(e_\b(z)\bigr)
\label{B0}
\ed
 for adjacent roots $\a$ and $\b$.
 {}From the definition of automorphisms and by simple use of \rf{3} we get
$$T_\b T_\a\bigl(e_\b(z)\bigr)\cdot K_\b^+\bigl(q^{2-c/2}z\bigr)^{-1}=$$
$$\frac{1}{(2\pi i)^2}\left[\oint_{C_{\infty}}\frac{du}{u}
\left(\oint_{C_{\infty}}\frac{dv}{v}e_\b(v)e_\a(q^2u)-
\oint_{C_{0}}\frac{dv}{v}\frac{q^{-1}v-q^2u}{v-qu}e_\a(q^2u)e_\b(v)\right)
f_\b(q^{2-c}z)-\right.
$$ $$
\left.\oint_{C_0}\frac{du}{u}f_\b(q^{2-c}z)
\left(\oint_{C_\infty}\frac{dv}{v}e_\b(v)e_\a(q^2u)-
\oint_{C_0}\frac{dv}{v}\frac{q^{-1}v-q^2u}{v-qu}e_\a(q^2u)e_\b(v)\right)
\right]=$$
 $$ \frac{1}{(2\pi i)^2(q-q^{-1})}\oint_{C_\infty}\frac{du}{u}
\left(\oint_{C_\infty}-\oint_{C_0}\right)
\frac{dv}{v}\frac{q^{-1}v-q^2u}{v-qu}\delta(v/q^2z)e_\a(q^2u)
K_\b^+(q^{2-c/2}z)-$$
\bn
\frac{1}{(2\pi i)^2(q-q^{-1})}\oint_{C_\infty}\frac{du}{u}
\left(\oint_{C_\infty}-\oint_{C_0}\right)
\frac{dv}{v}\delta(v/q^{2-2c}z)K_\b^-(vq^{c/2})e_\a(q^2u)
\label{B1}
\ed
Here the contour $C_{\infty}$ for the variable $u$ means that it enclose
 the origin and the point $u=qz$ is inside while for the contour $C_{0}$
 this point is outside. For the contour $C_{\infty}$ for the variable $v$
 the point $v=qu$ is inside while for $C_{0}$ it is inside.

The integration over $v$ cancel $\delta$ function and evaluate the integrand
 in a point $v=q^2z$ for the first integral and in a point $v=q^{2-2c}z$
 for the second. Then the integration over $u$ is a taking a residue in
 a point $u=qz$ for the first term (since the rest is normal ordered) and
 is zero for the second term.
As a result, we have
\bn
T_\a T_\b\bigl(e_\b(z)\bigr)=qe_\a(q^3z)
\label{B2}
\ed
which gives
\bn
T_\a T_\b T_\a\bigl(e_\b(z)\bigr)=qf_\a(q^{3-c}z)K^+_\a(q^{3-c/2}z)^{-1}
\label{B3}
\ed
Analogously,
\bn
T_\a T_\b\bigl(f_\b(z)\bigr)=qf_\a(q^3z)
\label{B4}
\ed
and
$$
T_\b T_\a T_\b\bigl(e_\b(z)\bigr)=T_\b T_\a\bigl(f_\b(q^{-cz})
K^+_\b(q^{-c/2}z)^{-1}\bigr)=
$$
\bn
qf_\a(q^{3-c}z)K^+_\a(q^{3-c/2}z)^{-1}.
\label{B5}
\ed
which proves \rf{B0}.

\setcounter{equation}{0}
\app{ Quadratic terms of $\RR $}

 We want to compare the quadratic terms of the expressions \rf{En14} and
 \rf{Pa21}. Let us for simplisity of notations restrict ourselves to
 degree $(0,0)$ terms and drop the common factor $(q-q^{-1})^{-2}$.

 As we noted before,  we can compare only those
 series in generators of the algebra, which are given in normal ordered form
 with respect to the grading index. The ordering
 can be given by the following  rules, which are equivalent to
 defining relations \rf{1}, \rf{2} (see \cite{KT5}):
$$e_{n+2k}e_{n}=q^2e_{n}e_{n+2k}+(q^{4}-1)\left(e_{n+1}e_{n+2k-1}+
q^{2} e_{n+2}e_{n+2k-2}+\ldots\right.$$
\bn
\left.+q^{2(k-2)} e_{n+k-1}e_{n+k+1}\right)+
(q^{2}-1)q^{2(k-1)}e_{n+k}^2 \qquad \mbox{for any}\ k>0,
\label{A2}
\ed
$$
e_{n+2k+1}e_{n}=q^2e_{n}e_{n+2k+1}+(q^{4}-1)\left(e_{n+1}e_{n+2k}+
q^{2} e_{n+2}e_{n+2k-1}+\ldots\right.$$
\bn
\left.+q^{2(k-1)} e_{n+k}e_{n+k+1}\right) \qquad\mbox{for any}\ k\geq 0,
\label{A3}
\ed
$$
f_{n+2k}f_{n}=q^{-2}f_{n}f_{n+2k}+(q^{-4}-1)\left(f_{n+1}f_{n+2k-1}+
q^{-2} f_{n+2}f_{n+2k-2}+\ldots\right.$$
\bn
\left.+q^{-2(k-2)} f_{n+k-1}f_{n+k+1}\right)
+(q^{-2}-1)q^{-2(k-1)}f_{n+k}^2 \qquad \mbox{for any}\ k>0,
\label{A4}
\ed
$$
f_{n+2k+1}f_{n}=q^{-2}f_{n}f_{n+2k+1}+(q^{-4}-1)\left(f_{n+1}f_{n+2k}+
q^{-2} f_{n+2}f_{n+2k-1}+\ldots\right.$$
\bn
\left.+q^{-2(k-1)} f_{n+k}f_{n+k+1}\right)\qquad\mbox{for any}\ k\geq 0.
\label{A5}
\ed
 The calculation of the integral
\bn
\oint_{C_1}\frac{dz_1}{z_1}\oint_{C_2}\frac{dz_2}{z_2}
  e(z_1)e(z_2)\ot f(z_1)f(z_2)
\label{A7}
\ed
 for the contours, enclosing the origin and inside  a region $|q^2z_2|< |z_1|
 < |q^{-2}z_2|$ can be performed as follows: we take first the integral in
 the region of analyticity $|z_1|>>|z_2|$ and substruct an integral over
 $z_2$ of a residue  (over $z_1$) in a point $z_1=q^{-2}z_2$.

 The regular term gives an expression
$$\sum_{n,m\in\ZZ} e_{-m}e_{-n}\ot f_mf_n.$$
 The corresponding normal ordered expression for its $(0,0)$ component
 has a form:
$$
 e_{0}^2\ot f_{0}^2+
(1-q^{-2})\sum_{k>0}q^{2k} e_{0}^2\ot f_{-k}f_{k}+
(1-q^2)\sum_{k>0}q^{-2k} e_{-k}e_{k}\ot f_{0}^2 +
(q^{2}+q^{-2})\sum_{k>0} e_{-k}e_{k}\ot f_{-k}f_{k}+
$$
\bn
(q^{2}-q^{-2})\sum_{k,l>0,l<k}q^{2l} e_{-k+l}e_{k-l}\ot f_{-k}f_{k}+
(q^{-2}-q^{2})\sum_{k,l>0,l<k}q^{-2l}  e_{-k}e_{k}\ot f_{-k+l}f_{k-l}.
\label{A6}
\ed
The residue over $z_1$ at a point $z_1=q^{-2}z_2$ comes from the simple pole
 of $f(z_1)f(z_2)$ at this point. So, we should take the residue
 of $f(z_1)f(z_2)$ and multiply by the evaluation $e(q^{-2}z_2)e(z_2)$.
The first expression can be given as a difference of two integrals:
$$
:f(q^{-2}z)f(z):=$$
$$
 \frac{1}{2\pi i}\oint_{C_1}f(z_1)f(z)\frac{dz_1}{z_1}
-\frac{1}{2\pi i}\oint_{C_2}\frac{z_1-q^2z_2}{q^2z_1-z_2}
f(z)f(z_1)\frac{dz_1}{z_1}
$$
where both contours enclose the origin and the point $z_1=q^{-2}z$ is inside
 the first contour and outside the second. Each of them can be derived by
integration of corresponding Laurent series. This gives
$$
:f(q^{-2}z)f(z):=
(q^{-2}-q^2)\sum_{n\in \ZZ}u^{-2n-1}q^{2n+2}
\left(\sum_{k\geq 0}q^{2k}f_{n-k}f_{n+k+1}\right)+
$$
\bn
(q^{-2}-q^2)\sum_{n\in \ZZ}u^{-2n}q^{2n}
\left((1-q^2)f_n^2+(q^{-2}-q^2)\sum_{k>0}q^{2k}f_{n-k}f_{n+k}\right).
\label{A8}
\ed
The evaluation of $e(q^{-2}z)e(z)$ assumes the direct application of
\rf{A2}, \rf{A3}. In an assumption $|q|<1$ it gives an analogous expression:
$$
e(q^{-2}z)e(z)=
\sum_{n\in \ZZ}u^{-2n-1}q^{2n+2}
\left(\sum_{k\geq 0}q^{-2k}e_{n-k}e_{n+k+1}\right)+
$$
\bn
\sum_{n\in \ZZ}u^{-2n}q^{2n}
\left(\frac{e_n^2}{1+q^2}+\sum_{k>0}q^{-2k}e_{n-k}e_{n+k}\right).
\label{A9}
\ed
For $|q|>1$ one should first perfom an analytical continuation (since
 in \rf{A9} we summed up the geometric progressions over $q$), that is
to use instead of $e(u)e(v)$ its analitycal continuation
$(u-q^2v)/(q^2u-v)e(v)e(u)$ and then do the same.
 The answer will coincide with \rf{A9}.
 The integral of the residue reduces to the  integration of a tensor
product of power
 series \rf{A8} and \rf{A9}. Its $(0,0)$- degree component has a form:
$$
\frac{(1-q^2)}{1+q^2} e_0^2 \ot f_0^2+ (1-q^2)\sum_{k>0}
q^{-2k}e_{-k}e_k \ot f_0^2+
$$
\bn
(q^{-2}-1)\sum_{k>0}q^{2k} e_0^2\ot f_{-k}f_k+
(q^{-2}-q^2)\sum_{k,p>0}q^{2(k-p)} e_{-p}e_{p}\ot f_{-k}f_k
\label{A10}
\ed
The half of the difference \rf{A6} and \rf{A10}
is
$$
{R}_2^0=\frac{1}{1+q^{-2}} e_{0}^2\ot f_{0}^2+
(1-q^{-2})\sum_{k>0}q^{2k} e_{0}^2\ot f_{-k}f_{k}+
q^{2}\sum_{k>0} e_{-k}e_{k}\ot f_{-k}f_{k}+
$$
\bn
(q^{2}-q^{-2})\sum_{k,l>0,l<k}q^{2l} e_{-k+l}e_{k-l}\ot f_{-k}f_{k}.
\label{A11}
\ed
 On the other hand, we can pick up the quadratic term of the infinite product
\rf{Pa18} modulo factor $(q-q^{-1})^{-2}$:
\bn
{R'}_2^0=\sum_{n,m\in \ZZ\atop n<m} e_{-n}e_{-m}\ot f_{n}f_{m}+
\frac{1}{1+q^{-2}}\sum_{n\in\ZZ} e_{-n}^2\ot f_{n}^2
\label{A1}
\ed
Performing normal ordering in the first tensor component, we get the
expression, which $(0,0)$ degree term  coincides with \rf{A11}.

\end{document}